\documentclass[]{gCFD2e}

\usepackage{subfigure}
\usepackage{algorithm, algorithmic}

\begin{document}



\title{Free Form Deformation, mesh morphing and reduced order methods: 
	enablers for efficient aerodynamic shape optimization}

\author{F. Salmoiraghi$^{\rm a}$, 
	A. Scardigli$^{\rm b \rm c}$$^{\ast}$
	\thanks{$^\ast$Corresponding author. Email: angela.scardigli@optimad.it \vspace{6pt}}, 
	H. Telib$^{\rm b}$ and G. Rozza$^{\rm a}$\\
	\vspace{6pt}  
	$^{a}${\em{SISSA, International School for Advanced Studies, mathLab, via Bonomea 265, 34136 Trieste, Italy}};
	$^{b}${\em{OPTIMAD Engineering srl, via Giacinto Collegno 18, 10143 Torino, Italy}};
	$^{c}${\em{Politecnico di Torino, DISMA -- Department of Mathematical Sciences, Corso Duca degli Abruzzi 24, 10129 Torino, Italy}};
\\\received{March 2018} }

\maketitle

\begin{abstract}
In this work we provide an integrated pipeline for the model order reduction of turbulent flows around parametrised geometries in aerodynamics.
In particular, Free-Form Deformation is applied for geometry parametrisation, whereas two different reduced-order models based on Proper Orthogonal Decomposition (POD) are employed in order to speed-up the full-order simulations: the first method exploits POD with interpolation, while the second one is based on domain decomposition. 
For the sampling of the parameter space, we adopt a Greedy strategy coupled with Constrained Centroidal Voronoi Tessellations, in order to guarantee a good compromise between space exploration and exploitation.

The proposed framework is tested on an industrially relevant application, i.e. the front-bumper morphing of the DrivAer car model, using the finite-volume method for the full-order resolution of the Reynolds-Averaged Navier-Stokes equations.


\end{abstract}

\section{Introduction}\label{sec:Intro}
We would like to start from a simple calculation to give the motivation of the present work.
Let us take a general, fluid dynamic industrial problem.
Typically, the Reynolds number requires to exploit turbulence models to solve these kind of problems.
The industrial standard is nowadays the Reynolds Averaged Navier-Stokes (RANS) closure model.
This algorithm may have a cost of $O$(10$^2$--10$^4$) cpu hours for a standard 3D simulation employing a computational grid of $O$(10$^7$--10$^8$) elements.
Even if the license cost is null, the energy consumption cost of modern computers is about 0.05 euro per cpuh, that is, $O$(5--5$\times10^2$) euro per each simulation.
In an evolutionary optimization campaign, for non trivial problems, we perform $O$(10$^2$--10$^3$) simulations, that leads to a cost of $O$(5$\times10^2$--5$\times10^5$) euro .
Automatic shape optimization is used only if strategic and addressed in a reasonable amount of time.

Computer performance is constantly improving and the cost of simulations should, theoretically, decrease.
In practice, however, it can be observed that numerical models always tend to saturate computational resources, since more complex and accurate models are continuously developed.
Computational time does not decrease as computing power increases, and this represents a big challenge for optimization.

We can rely on the observation of physical phenomena to develop new methodologies to overcome the mentioned limitations.
When a physical problem depends on some changing physical/geometrical parameters, it often happens that the solution changes smoothly.
Starting every time from scratch to solve the problem is not optimal or even feasible.
The intuition is to use a few high-fidelity simulations for some properly selected values of the parameters (related to different configurations) to build a solutions database and then to recycle problem data, solving it by combining the solutions in the database.
This is nothing but the rationale that drives the development of every Reduced Order Method (ROM).
Thanks to ROMs, we are able to perform simulations of complex phenomena almost in real time, after the construction of the database containing $O(10)$ high-fidelity, indeed very expensive, solutions of the problem for properly selected parameter values.

An overview of different reduced order methods strategies for industrial and biomedical applications can be found in~\cite{salmoiraghi2016advances} and \cite{bergmann2014reduced}.

In the following we present a new pipeline and different tools for model reduction of complex and industrial problems.
In Section~\ref{sec:Geometry Morphing} we present two different approaches to the geometry parametrization, exploiting Free Form Deformation (FFD), then, after Section~\ref{sec:FOM} devoted to full order methods, in Section~\ref{sec:ROM} we present two different reduced order methods. 
In both Sections, we give a theoretical insight and then some algorithm details for the introduced methodology, as well as a properties comparison between them.
Finally, in Section~\ref{sec:Results}, we apply the presented methodologies to an automotive industrial benchmark, namely the \emph{DrivAer}~\footnote{http://www.drivaer.com} model.
This model presents interesting features both from the physical and geometrical point of view.
In fact, on the one hand, the high Reynolds number implies the use of a RANS model, on the other hand, the geometry is also made up of rotating wheels that can not be neglected, in order to obtain accurate results, but this increases the complexity of the phenomena involved.
Both the turbulence models and rotating mechanical parts represent quite new challenges in the ROM community.
This field is quite new in the reduced order methods community and this is one of the aspect of novelty of the present work.

\section{Geometry Morphing}\label{sec:Geometry Morphing}
In the present work, we provide an approach to perform geometry parametrization and deformation based on the Free Form Deformation (FFD) method, as described in its original version in~\cite{sederberg1986free}.
Starting from this pioneering work, which is related to solid geometric modelling, FFD has been developed mainly in the frame of computer graphics and only recently employed in aerodynamic shape design problems.

\begin{figure}
	\centering
	\includegraphics[width=0.6\textwidth]{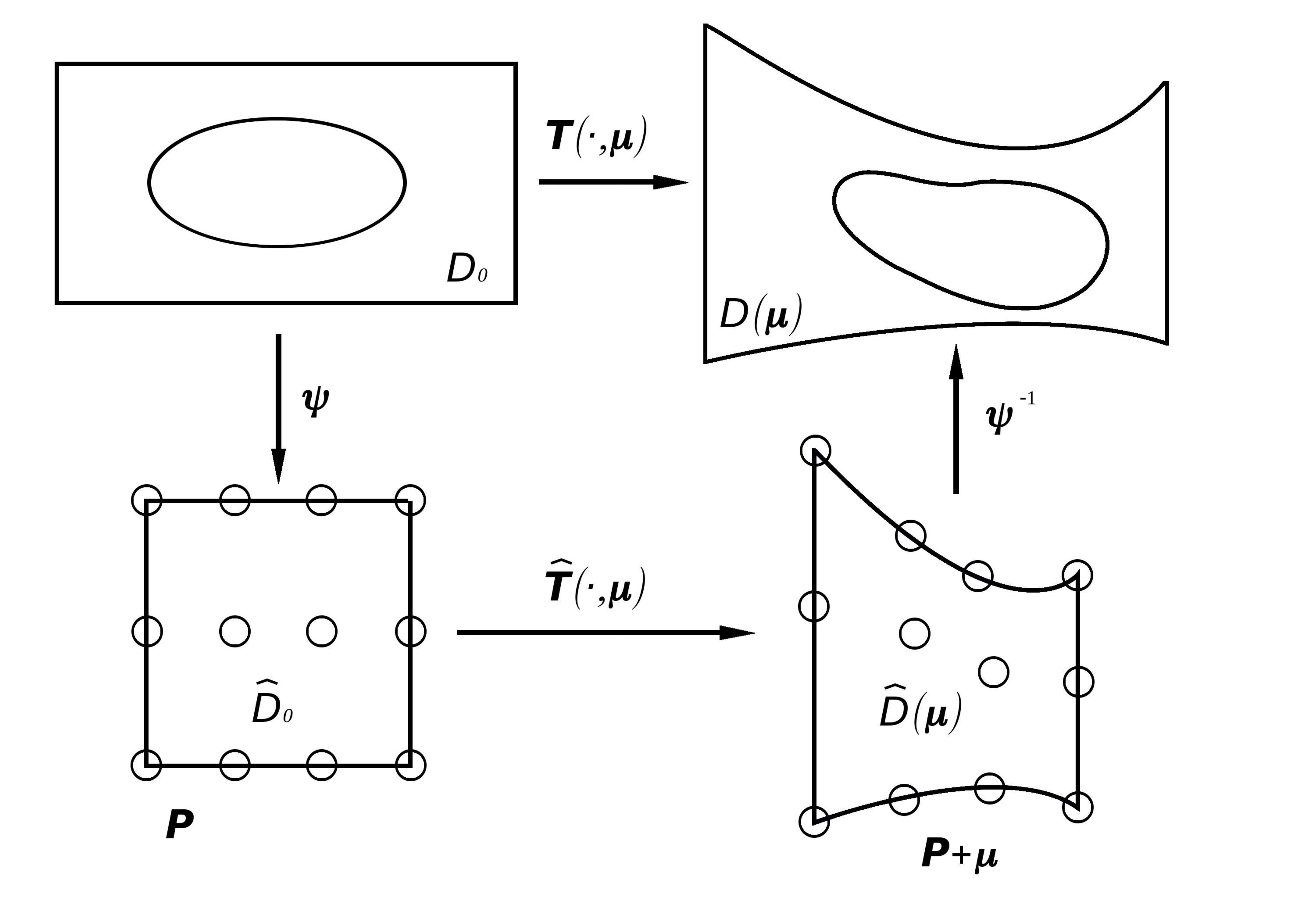}
	\caption{Sketch of the FFD map construction.}
	\label{fig:FFD sketch}
\end{figure}

Basically, this technique first sets a control lattice surrounding the part of the geometry to be morphed and then deforms the geometry in a continuous and smooth way by moving only the control points of such lattice.
This operation is carried out in three steps, shown in Figure~\ref{fig:FFD sketch}: first we need to map the actual domain $D_0$ to the reference one $\hat{D}_0 = [0,1]^d$ (where $d$ is the physical dimension of $D_0$) through the map $\bm{\psi}$; second we set a regular grid of unperturbed control points $\bm{P}$ and we move some control points of a quantity $\bm{\mu}$ to deform the lattice through the map $\hat{\bm{T}}$; finally, we map back the resulting domain $\hat{D}$ to the physical domain $D$ through the inverse map $\bm{\psi}^{-1}$.
The FFD map is the combination of the three maps:
\begin{equation}
\bm{T}(\cdot, \bm{\mu}) = (\bm{\psi}^{-1} \circ \hat{\bm{T}} \circ \bm{\psi}) (\cdot, \bm{\mu}).
\end{equation}
For a deeper insight, see \cite{lassila2010parametric, koshakji2013free, forti2014efficient} 
and \cite{sieger2015shape} for an application very similar to ours. A comprehensive dissertation on different geometry deformation techniques (in particular FFD based techniques) can be found in~\cite{anderson2012parametric}. 

In the present work, we apply the FFD paradigm to morph directly the computational mesh around the geometrical model.
In the following we provide some insight about the general method, whereas in Section \ref{subsec:PyGeM} we focus on its application to our problem, highlighting advantages and drawbacks.

\subsection{Morphing features}\label{subsec:features}
Let us now provide some general considerations about the FFD technique.
Since FFD is independent from the topology of the object to be morphed, it is extremely versatile and suitable to parametrize very complex geometries, including volume meshes, surface triangulations and CAD representations.
Moreover it permits to obtain large deformations as well as small ones.
Such strength of the method can turn into a weakness when handling 2D or 3D meshes: the input mesh should be good enough, in terms of number and quality of mesh elements, to guarantee a good description of the deformed geometry. 
Low quality inputs will results inevitably in poor deformations.
This behaviour can be overcome, for instance, if we perform FFD on the manifold of the geometry, e.g for the morphing of CAD geometries or as in~\cite{salmoiraghi2016isogeometric} for the morphing of Isogeometric geometries (NURBS surfaces); alternatively, some curvature based refinement needs to be introduced for mesh adaptation. 

\begin{figure}
	\centering
	\includegraphics[height=0.3\textheight]{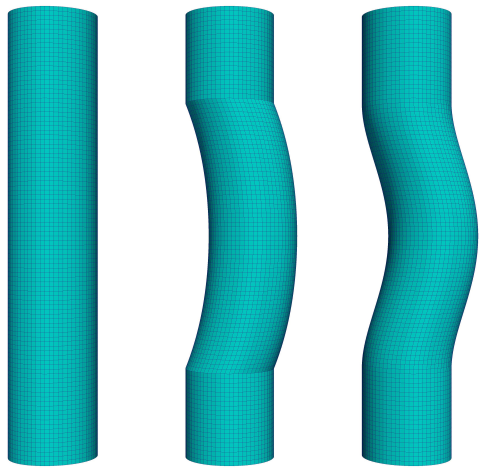} 
	\caption{Regularity of the deformation: original mesh (left), $\mathbb{C}_0$ deformation (center) and $\mathbb{C}_1$ deformation (right).}
	\label{fig:C0 C1 mesh}
\end{figure}

As it is, the method is suited for global shape deformations and when the geometry is unconstrained. Generally speaking, this is not the case for real-life application, where local deformations may be involved as well as the need to impose arbitrary shape constraints.
This means that we need some control over the continuity and smoothness of the deformation: even if it is driven by the properties of Bernstein polynomials inside the FFD control box, we can encounter some problems in the interface between the deformed and undeformed part of the domain.
To overcome this problem, we can insert some control points to be kept fixed close to the interface.
For the sake of clarity, we provide a simple example.
Consider the mesh on a cylinder in Figure~\ref{fig:C0 C1 mesh} (a). If we construct the FFD control lattice with three points in the vertical direction and we move only the second one, the results is depicted in Figure~\ref{fig:C0 C1 mesh} (b). 
The mesh in this case has two artificial edges and is $\mathbb{C}_0$ continuous.
On the other hand, if we use five control points in the vertical direction and move only the third one, we obtain a smoother, indeed $\mathbb{C}_1$, mesh.
A more sophisticated approach (~\cite{scardigli2015enabling}) consists in introducing a filter scalar function (with the required order of continuity) that weights the deformation according to the distance from the constraint itself.

Another feature that comes for free with FFD is the straightforward parallelization of the code. 
In fact, we can ideally move independently each point of our geometrical object on a different CPU, passing only a limited amount of information among processors.
Nevertheless, this may not be the case for the mesh morphing of large-scale problems, where the implementation of deformation controls is more complex. 

Thanks to these features, the FFD is a valid option for industrial problems: however, the choice of the preferable FFD approach should depend on the actual application, as pointed out in Section \ref{subsec:PyGeM}. 

\subsection{Mesh morphing strategy}\label{subsec:PyGeM}
As stated before, we perform FFD directly on a computational mesh around the model.
In such a way, provided the mesh for the reference geometry, we extract and map (through FFD) only the coordinates of the vertices, leaving the connectivity and the other properties of the mesh untouched.

In view of reduction strategies, this framework allows also to obtain the mapping of the solution for free, so that we have the value of the solution on the same degrees of freedom without introducing an interpolation of the solution between different meshes.

In terms of deformation control, the $\mathbb{C}_0$ continuity between the fixed and the deformable part of the mesh can be easily fulfilled if we do not move the external control points of the FFD bounding box.
The non-penetration condition of the cells, instead, is guaranteed by limiting the displacements in order to avoid the overlapping of the control points.  
However, depending on the deformation type, different drawbacks may affect the cells (e.g. high skewness, high non-orthogonality, etc. ) and impair the quality of the simulation results.
To prevent such behaviour, a set of application-specific mesh quality constraints has to be satisfied: roughly speaking, the more problematic cells tends to be the ones close to the interface between the deformed and undeformed part of the domain.
Increasing the continuity of the deformation allows to overcome this potential problem.
This can be done by adding more control points to be kept fixed close to the interface, as explained in Section~\ref{subsec:features}.
Despite these expedients, guaranteeing the minimum quality of the volume mesh may be infeasible or really difficult to achieve in those cases where non-small deformations are required or the parametrization of rather complex geometries is involved.

In the case of shape optimization, aimed at finding the optimal configuration for small deformation of the starting geometry, this drawback is not a great limitation.
Moreover, FFD on mesh allows to skip the mesh generation for every new configuration, leading to  significant time savings.
In fact, even if on large meshes the deformation itself and the quality checks may be time consuming, the mesh generation is usually more demanding in terms of computational time and resources.

This strategy relies on PyGeM~\footnote{Python {G}eometrical {M}orpher -- https://github.com/mathLab/PyGeM/}, which is a \textit{Python} library using Free Form Deformation and Radial Basis Functions to parametrize and morph complex geometries, including CAD representations. 

\section{Full Order Model}\label{sec:FOM}
In the present work, the full-order model is represented by the resolution of Navier-Stokes equations (NSE). 
In broad terms, the NSE are a system of time-dependent, non-linear, partial differential equations which govern the motion of fluids.
Further difficulties arise when turbulence is involved, as occurs in many engineering applications: turbulent flows exhibit a chaotic behaviour, characterised by significant and irregular variations in space and time, and their study represents a challenge under both the analytical and numerical point of view (see for instance \cite{pope2011turbulent}, for a more comprehensive review of the problem). 

Nowadays, there are four main approaches to deal with turbulence: Direct Numerical Simulation (DNS), Reynolds-Averaged Navier-Stokes (RANS), Large Eddy Simulation (LES) and hybrid LES-RANS models.
In DNS, all the scale of motions are solved for one realization of the flow.
Although easy in principle, solving the whole range of spatial and temporal turbulence scales is often not feasible, given the complexity of the phenomena. 
DNS is indeed very expensive and the computational costs tends to increase cubically with the Reynolds number $Re$: therefore, this approach can be applied to flows characterized by low or moderate $Re$, whereas it has prohibitive costs for industrial applications at higher $Re$.
Other techniques for simulating turbulent flows, such as LES (\cite{sagaut2006large,pope2011turbulent}) and hybrid models (\cite{frohlich2008}), have started to be employed in engineering applications. 
Nevertheless, the resolution of RANS equations is still the most common approach in industry, especially in early stages of design or during aerodynamic optimization, when several simulations are required.
The general idea behind the RANS approach is to decompose velocity $\bm{U}$ and pressure $p$ into ensemble-averaged and fluctuating components (Reynolds decomposition), obtaining approximate solutions to the NS equations. 
A turbulence model is required to determine the Reynolds stress, the unknown term which accounts for fluctuations contribution, and to provide closure to the system of equations (\cite{pope2011turbulent}).

In the following, we will refer to the RANS equations as the high-fidelity/full-order model.

\section{Reduced Order Model}\label{sec:ROM}
We proposed two different approaches for the construction of the reduced basis, both relying on the Proper Orthogonal Decomposition (POD).
The POD allows to extract the modes from a set of solutions of the problem at hand for different values of the parameters. 
For a deeper insight on POD techniques, see~\cite{Ravindran99properorthogonal, POD_aubry_1991, chinesta2011short, chinesta2016model, chinesta2014separated} for the general formulation; see~\cite{benner2015survey} for dynamical systems. For the application of POD to the model order reduction of fluid dynamics problems, we recall~\cite{manzoni2015reduced} for potential flows, \cite{salmoiraghi2016isogeometric, ballarin2015supremizer} for viscous flows and~\cite{bui2003proper} for compressible flows.

\subsection{POD at a glance}\label{subsec:POD}
The idea is to start from a parametric geometrical model and create a database $\bm{\Xi} = [\bm{\mu}_1 \ | \  \ldots \ | \ \bm{\mu}_N]$ of parameter values and a database $\bm{\Theta} = [ u(\bm{\mu}_1) \ | \  \ldots \ | \ u(\bm{\mu}_N)]$ of outputs thanks to a proper sampling strategy, shown in section~\ref{subsec:Sampling}.
Once we have the database we perform the Singular Value Decomposition (SVD) of the sample to extract the POD modes $\bm{\psi}$
\begin{equation}\label{eq:svd}
\bm{\Theta} = \bm{\Psi} \bm{\Sigma} \bm{\Phi}^T,
\end{equation}
where $\bm{\Psi}$ and $\bm{\Phi}$ are the left and right singular vectors matrices of $\bm{\Theta}$, and $\bm{\Sigma}$ is the diagonal matrix containing the singular values in decreasing order. The POD modes $\bm{\psi}$ are nothing but the columns of $\bm{\Psi}$.

Alternatively, it is possible to compute the basis through the Method of snapshots \cite{sirovich1987turbulence}, by solving the equivalent eigenvalue problem (\cite{volkwein2013proper})
\begin{equation}\label{eq:sirovich}
\bm{\Theta}^T\bm{\Theta} \phi_i = \lambda_i \phi_i,
\end{equation}
and setting
\begin{equation}
\psi_i=\frac{1}{\sqrt{\lambda_i}}\bm{\Theta} \phi_i,
\end{equation}
When the snapshot dimension is much greater than $N$, it is less expensive to solve Problem~\ref{eq:sirovich}, whereas the first approach is more reliable for badly conditioned matrices (\cite{demmel1997applied}).

Exploiting the new basis, we can express the reduced solution of the problem as:
\begin{equation}
u^N = \sum_{i=0}^N \alpha_i \psi_i,
\end{equation}
that is, a combination of the basis functions. The problem shift to find, for each new value of the parameter, the value of the coefficients $\alpha_i$. 
Some (\cite{rozza2008reduced, quarteroni2011certified, hesthaven2015certified}) use the weak formulation of the problem and exploit the POD modes as basis functions to find the coefficients (suitable for problems with weak formulation and affine dependence from the parameters), others (\cite{peherstorfer2015online, peherstorfer2015dynamic}) use the measurement of the quantity of interest coming from sensors (suitable for real time evaluation `on the field').
In the following we show two different methods for the POD coefficients evaluation exploiting two different paradigms: Proper Orthogonal Decomposition with Interpolation (PODI) and Proper Orthogonal Decomposition combined with Domain Decomposition techniques (DD-POD).

\subsection{PODI}\label{subsec:Scenario 1}
\begin{figure}
	\centering
	\includegraphics[scale=.06]{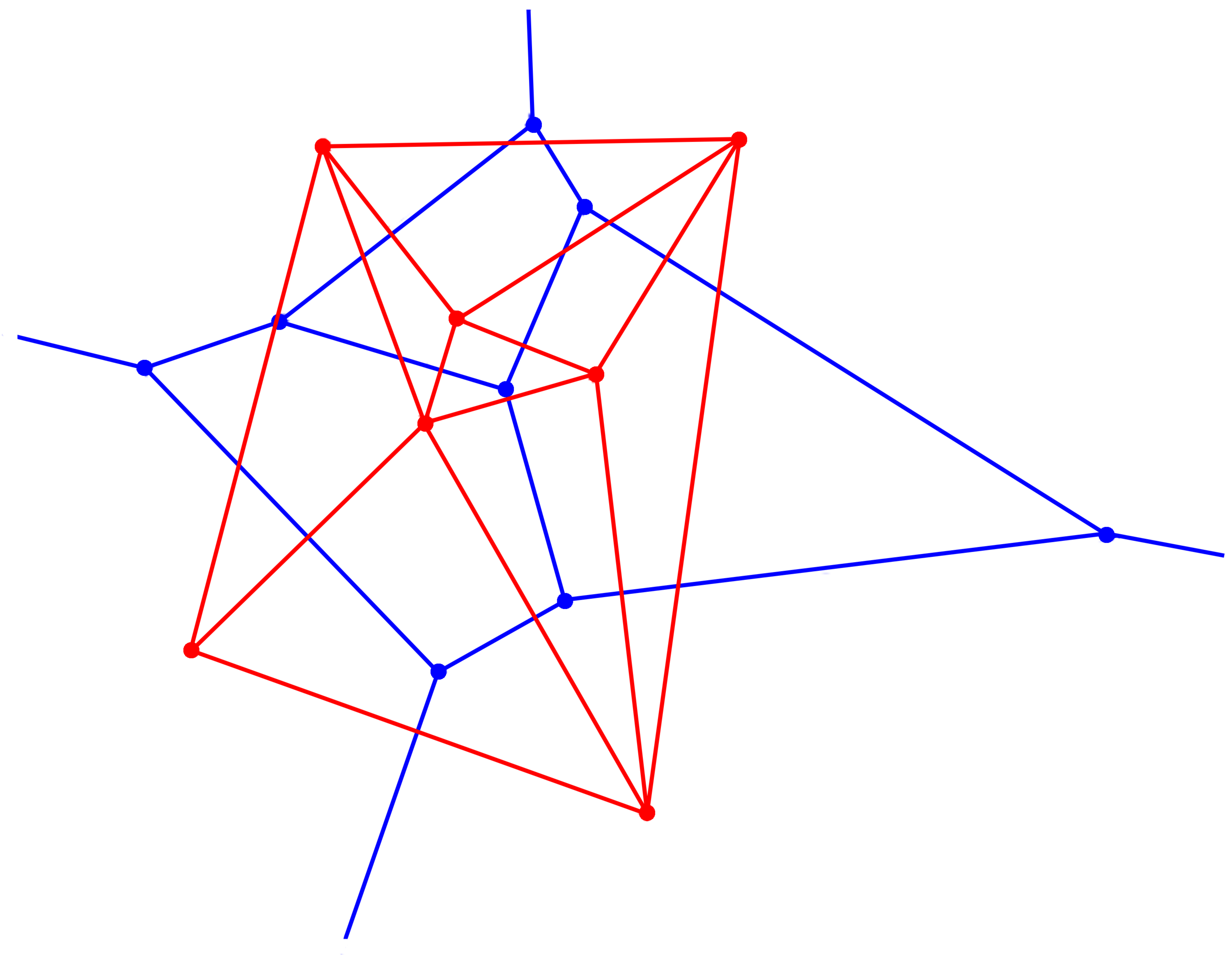}
	\caption{The duality between Voronoi tesselation (blue) and Delaunay triangulation (red).}
	\label{fig:Tesselation}
\end{figure}

PODI was first introduce in~\cite{tan2003proper} and it has been used recently in aerodynamic applications (\cite{dolci2016proper}).
The rationale behind PODI regards the evaluation of the POD coefficients by interpolation of the POD coefficients computed for the parameter points $\bm{\mu}_k \in \bm{\Xi}$.
In these points, the reduced and high fidelity solutions ($u^N$ and $u$ respectively) are equal by construction: 
\begin{equation}
\forall \bm{\mu}_k \in \bm{\Xi}: \displaystyle u(\bm{\mu}_k) = u^N(\bm{\mu_k}) = \sum_{i=0}^N \alpha_i (\bm{\mu_k}) \psi_i
\end{equation}
For each new value of the parameter $\bm{\mu}_{new}$, we interpolate the $\alpha_i (\bm{\mu_k})$ coefficients to find the new $\alpha_i (\bm{\mu}_{new})$ coefficients and multiply them for the reduced basis matrix $\bm{\psi}$ in order to evaluate the new reduced solution
\begin{equation}
u_{new}^N = \sum_{i=0}^N \alpha_i (\bm{\mu}_{new}) \psi_i
\end{equation}
The interpolation is performed exploiting the \textit{Delaunay} triangulation, and its dual \textit{Voronoi} tessellation (see figure~\ref{fig:Tesselation}), of the parameter space (\cite{fortune1992voronoi, cammi2016reduced, du1999centroidal, watson1981computing}).
As a consequence, the method can be used efficiently in the N-th dimensional case.

This first model order reduction strategy relies on EZyRB~\footnote{Easy Reduced Basis method -- https://github.com/mathLab/EZyRB/} tool.
EZyRB is a \textit{Python} library for model order reduction based on barycentric triangulation for the selection of the parameter points (see Section \ref{subsec:Sampling} for algorithm details) and on POD for the selection of the modes. 
The term easy (EZy) is used with respect to the classical Reduced Basis Method: in RBM we reconstruct the whole solution of the problem thanks to its weak formulation and an hypothesis of affine dependence from the parameters; in EZyRB we reconstruct only the output of interest thanks to an interpolation of the reduced basis (POD modes) coefficients.
It is ideally suited for actual industrial problems, since its structure can interact with several simulation software simply providing the output file of the simulations.
This tool is suited for the construction of the model order reduction strategy when we start our sampling from scratch (since it is based on sampling strategy shown in Section~\ref{subsec:Sampling}), but it can be used in a sub-optimal way even if the database has been already computed in the past.

\subsection{DD-POD}
\label{subsec:Scenario 2}
\begin{figure}
	\centering
	\includegraphics[width=0.6\textwidth]{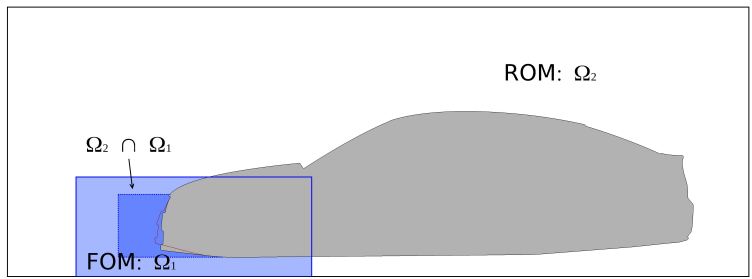}
	\caption{DD-POD example: $\Omega_1$ (blue), $\Omega_2$ (white) and overlapping region $\Omega_1\cap\Omega_2$ (light blue).}
	\label{fig:podFOAM framework}
\end{figure}

The main idea behind the method, as first proposed by~\cite{buffoni2009iterative}, is to split the domain in two parts and use a different approximation method in each region.
Generally speaking, POD based ROMs suffer from the fact that the basis space is a linear combination of the solution space spanned by the snapshots: this means that it is not capable to represent non-linear phenomena, unless a sufficient rich database is provided.
To bypass such difficulty, a hybrid low-order/high-order method based on domain decomposition is employed.
As shown in Figure~\ref{fig:podFOAM framework}, the computational domain is decomposed into two regions, i.e. $\Omega_1$ and $\Omega_2$ such that $\Omega_1 \cap \Omega_2 \neq \emptyset$: the canonical CFD solver is used where the effects of non-linearities and geometry variations are predominant ($\Omega_1$), whereas linear and weakly non-linear phenomenology is addressed by the ROM ($\Omega_2$).
The two models are coupled through the overlapping region using a Schwarz-type method, resulting in a non-local boundary condition on $\partial \Omega_1$. 
The main steps of the method are presented in Algorithm~\ref{alg:PODwDD}.
\begin{algorithm}
	\caption{DD-POD algorithm}\label{alg:PODwDD}
	\begin{algorithmic}[1]
		\STATE set initial boundary conditions for $u_1(\bm{\mu}_{new})$ on $\partial \Omega_1$
		\WHILE{$\textit{convergence} =\text{false}$}{}
		\STATE evaluate $u_1$ by integrating the governing equations in $\Omega_1$
		\STATE $\alpha \gets \arg\,\min_\alpha \left( \parallel u_1-\sum_{i=0}^N \alpha_i \psi_i\parallel_{\Omega_1\cap\Omega_2}\right)$
		\STATE $u_2(\bm{\mu}_{new})^N \gets \alpha$
		\STATE evaluate $u_2^N|_{\partial \Omega_1}$
		\STATE update boundary conditions for $u_1$ on $\partial \Omega_1$ 
		\STATE check for \textit{convergence}
		\ENDWHILE
	\end{algorithmic}
\end{algorithm}

At every solver iteration, the high-order solution $u_1$ in $\Omega_1$ is used to evaluate the POD coefficients $\alpha_i$, by solving a least-square problem which minimizes the L2-norm of the distance between the CFD and ROM solutions in $\Omega_1\cap\Omega_2$.
This allows to reconstruct the POD solution $u_2^N$ in $\Omega_2$, and in particular its restriction to $\partial \Omega_1$, determining a new set of boundary conditions for the CFD solver.

Such operation has a negligible cost with respect to the iteration and can be integrated in the CFD solver with little effort: for the current application, the algorithm is implemented directly in the \textit{OpenFOAM\textsuperscript{\textregistered}} solver.

\subsection{Sampling strategy}
\label{subsec:Sampling}

By construction, the POD basis gives an optimal representation, in terms of energy, of the solution space. 
This means that the error of both ROMs is ideally zero for each snapshot belonging to the database: in the limit of no compression, the basis is capable to reproduce exactly the solution.
However, a robust ROM should ensure a sufficient accuracy over the entire parameter space, i.e. for each $\bm{\mu}_{new}$: to fulfil this requirement, a potentially large number of snapshots may be necessary.
Since the computational cost of such evaluations is often prohibitive for industrial applications, the use of an efficient sampling strategy on the parameter space assumes high relevance in this context.

In the present work, we adopt the approach first proposed by~\cite{lombardi2011low}, which couples Constrained Centroidal Voronoi Tessellations (CCVT) and Greedy methods.
In this strategy, new well-spaced points are added iteratively, enriching the database in those regions  where a certain error indicator exceeds a fixed tolerance.
Our error indicator, i.e. the density function used in the CCVT, is built exploiting a Leave-One-Out Cross-Validation technique.
Given an initial database containing the snapshots for the vertices of the parametric domain $\bm{\Theta}_0 = [ u(\bm{\mu}_1) \ | \  \ldots \ | \ u(\bm{\mu}_{N_0})]$, we compute $N_0$ POD-bases, one for each $u(\bm{\mu}_k)$ in $\bm{\Theta}_0$, by leaving that snapshot out: it is then possible to evaluate $u^N(\bm{\mu}_k)$ as the projection of $u(\bm{\mu}_k)$ on the POD-basis spanned by all the remaining snapshots and compute our indicator $e_s = \parallel u(\bm{\mu}_k) - u^N(\bm{\mu}_k) \parallel$.
Such information is used to compute the centroids of the tessellation elements (or the barycentric values of the Delaunay triangulation, which is its dual): among those points, we choose the one where the density function reaches its highest value as the new sampling point. 
The strategy is summarised in Algorithm~\ref{alg:leave1out}.

\begin{algorithm}
	\caption{Sampling strategy with leave-one-out algorithm}\label{alg:leave1out}
	\begin{algorithmic}[1]
		\STATE compute $\bm{\Theta}_0$
		\WHILE{$\max e_s > tol$}{}
		\FORALL{$\bm{\mu}_k$ in $\bm{\Xi}$ }
		\STATE remove $u(\bm{\mu}_k)$ from $\bm{\Theta}$
		\STATE compute ROM
		\STATE evaluate  $u^N(\bm{\mu}_k)$
		\STATE error $e_s(k) = \parallel u(\bm{\mu}_k) - u^N(\bm{\mu}_i)^N \parallel$
		\ENDFOR
		\FORALL{\textit{simplex} in Delaunay triangulation}
		\STATE error $e_t = \texttt{Area} * \sum_{\text{vertices}} e_s$
		\ENDFOR
		\STATE refined simplex $\gets \arg\max e_t$ 
		\STATE $\bm{\mu}_{new} = \frac{\sum \bm{X} e_s}{\sum e_s}$
		\ENDWHILE
	\end{algorithmic}
\end{algorithm}

An important drawback of the methodology is represented by the curse of dimensionality.
In fact, in order to start the sampling algorithm, we need to compute the solution of all the vertices of the parametric domain, that is, $2^n$ solutions, where $n$ is the number of parameters involved.
This means that the computational cost for the evaluation of $\bm{\Theta}_0$ will grow exponentially with $n$, and may become prohibitive for large-scale industrial problems.
This characteristic is a bottleneck of the present framework.

\subsection{Critical comparison between the two methods}\label{subsec: ROM compare}
Let us now provide some considerations and comparisons between the two techniques.
Thanks to the simplicity of the method, with PODI we can build the reduced order model on the solution only in the region of interest (e.g. the surface of the body) and not in the whole computational domain. 
Therefore, PODI on output is faster and quicker both in the construction and evaluation of the reduced order model. 
In fact, during the construction step, we import, assemble and perform SVD on matrices containing the output that commonly are $O(10^-3)$ smaller than the ones containing the whole solution.
During the evaluation step, we only evaluate the new coefficients by interpolation and then we perform a matrix-vector product, which is a trivial operation.
On the other hand, with the second approach, we can reconstruct the whole, even non-linear, solution, since it partially solves the high-fidelity problem during the online phase.
If the field reconstruction is not sufficiently good, instead of adding new snapshots, it is possible to enlarge the domain where we solve the high-fidelity model and obtain a more accurate result: obviously, the more we extend the inner domain, the less we gain in terms of computation speed-up.
For such reasons, the second strategy typically requires fewer snapshots.
Moreover, the DD-POD approach does not rely on the parametrization: given a set of high-fidelity simulations, it can be used straightforwardly even if the snapshots are generated with different parametrizations.
This is not the case for the PODI approach, where the geometry parametrization has to be known in order to perform interpolation.
These features give a good flexibility to the second method, but the model reduction is less severe, leading to an online evaluation step slower than the one related to the first strategy.

\section{Results on DrivAer model}\label{sec:Results}
The results of this integrated approach, from geometrical morphing to model reduction are performed on the \textit{DrivAer} model, which is a realistic generic car model developed by TU Munich in collaboration with Audi AG and BMW Group, and made available in several configurations for research purposes.
Based on two medium-size cars (see~\cite{heft2011investigation, heft2012introduction} for further details), the model represents a good compromise between complex production cars and strongly simplified models.
Generic models like the Ahmed body (\cite{ahmed1984some}) or the SAE model (\cite{cogotti1998parametric}), are widely used in vehicle aerodynamics to investigate basic flow structures but fail to predict more complex phenomena, due to the oversimplification of the geometries in relevant regions, e.g. rear end, underbody and wheelhouses.
On the other hand, real cars are unlikely to be used for validation purposes, due to data access restrictions.
In this scenario, the DrivAer model constitutes a solid benchmark for industrial applications, merging a realistic and detailed geometry description with the availability of numerical and experimental validation data.

We adopt a fastback configuration with mirrors, rotating wheels and smooth underbody.
All the numerical investigations are carried out with ground simulation and at realistic Reynolds numbers, increasing the complexity of the phenomena:  the flow is fully three-dimensional and turbulent, characterized by separation, recirculation and unsteady wakes (\cite{hucho2013aerodynamics}).

\subsection{High-fidelity model: simpleFoam}\label{subsec:simpleFoam}
High fidelity evaluations are carried out using \textit{OpenFOAM\textsuperscript{\textregistered}}~\footnote{http://www.openfoam.org}, an open source software for computational fluid dynamics. More specifically, we decide to use \textit{simpleFoam}, the steady-state solver for incompressible flows based on a semi-implicit method for pressure-linked equations (SIMPLE) algorithm, implementing a cell-centred finite volume method for Reynolds-Averaged Navier-Stokes (RANS) equations (see~\cite{ferziger2002computational} for a comprehensive overview of solution methods).
It should be noted that the proposed approaches do not rely on the choice of the solution methods and software, but can be coupled with any canonical CFD solver. 

A three-dimensional hex-dominant mesh of about $15\cdot10^6$ cells is generated around the car model by the \textit{snappyHexMesh} utility, introducing symmetry in the longitudinal plane. 

The Reynolds number of the simulations, based on the free-stream velocity and on the car length, is set to $4.87\cdot10^6$.
The Spalart-Allmaras turbulence model (\cite{spalart1992oneequation}) is employed and we introduces wall-functions to describe the near-wall flow. 
Despite such expedients, the computational cost of each high fidelity evaluation is of $O(10^2)$ cpu hours on a last generation supercomputer, without considering the mesh generation.

Once we get the solution, we post-process it in order to extract the output of interest, namely the pressure on the surface $p_\mathrm{w}$, the wall shear stress $\bm{\tau}_\mathrm{w}$ and the drag coefficient, defined as
\begin{equation}
C_\mathrm{x} = \frac{\int _S -p \bm{n} \cdot \bm{x} + \int _S \bm{\tau} \cdot \bm{x}}{\frac{1}{2} \rho U^2 A_\mathrm{f}} 
\label{eq:force coeffs}
\end{equation}
where $\bm{n}$ is the outward-pointing versor normal to the model surface, $\bm{x}$ and $\bm{z}$ are the versors in the longitudinal and vertical directions, respectively, $\rho$ is the air density, $U$ the unperturbed speed and $A_\mathrm{f}$ is the frontal area of the model.

Compared with experimental data (see \cite{heft2012introduction}), our setup leads to a $5\%$ error on the output of interest, which is considered acceptable for the work purpose (and its applications).

\subsection{Geometry Morphing}\label{subsec:Results mesh morphing}
\begin{figure}
	\centering
	\includegraphics[width=.7\textwidth]{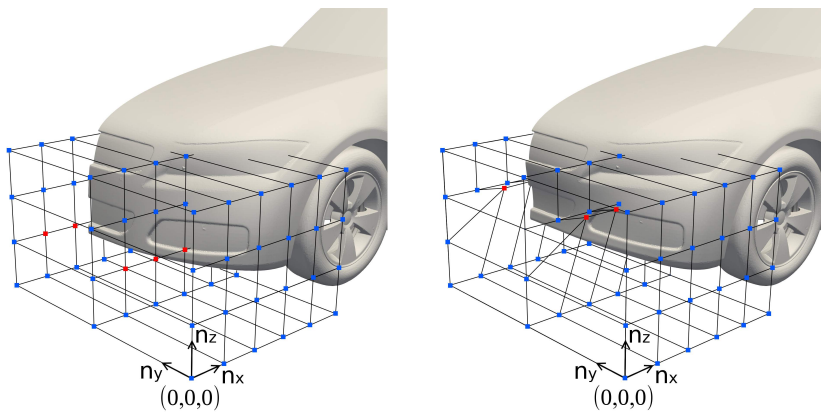}
	\caption{FFD lattice and a possible deformation (red control points): $\left[\mu_1,\mu_2 \right] = \left[0.18, 0.30\right]$.}
	\label{fig:Comparison with lattices}
\end{figure}

For the problem at hand, we wrap the mesh close to the car front bumper in a lattice of $6\times3\times4$ control points and we keep fixed all the points except the ones highlighted in red in Figure~\ref{fig:Comparison with lattices} (points with indices $[n_x, n_y, n_z] = [1,1,1] \lor [1,2,1] \lor [2,1,1] \lor [2,2,1] \lor [3,1,1] \lor [3,2,1]$).
Such points are allowed to move by the same quantity in the longitudinal and vertical directions, resulting in a 2-dimensional parameter space, $\left[\mu_1, \mu_2 \right]$, with bounds  $\left[\mu_1, \mu_2 \right]_{min} = [-0.18, -0.3]$ and $\left[\mu_1, \mu_2 \right]_{max} = [0.18, 0.3]$, i.e. a little less of the distance between two consecutive points in the corresponding direction.
The chosen parametrization guarantees the overall satisfaction of our set of mesh quality constraints.

\begin{figure}
	\centering
	\includegraphics[width=0.5\textwidth]{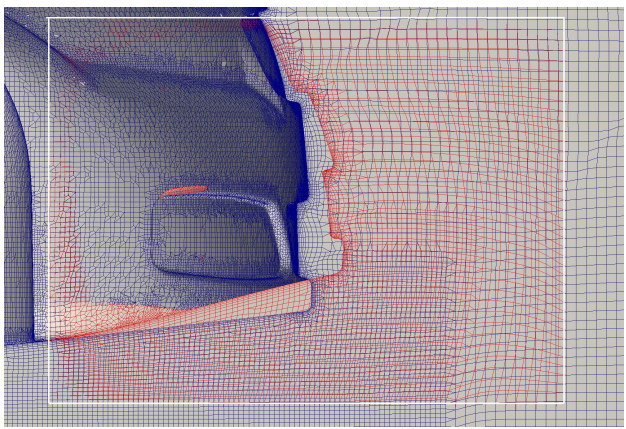}
	\caption{Mesh morphing of the DrivAer front bumper: original mesh (blue lines) vs modified mesh (red lines) for $\left[\mu_1, \mu_2 \right] = [-0.18, 0.3]$. The white box identifies the deformed region.}
	\label{fig:meshmorphing}
\end{figure}

In Figure~\ref{fig:meshmorphing}, we show the resulting deformed mesh for one of the vertices of the parameter space.
Since the FFD is applied directly to the volumetric mesh, we do not need to rebuild the mesh for the others configurations, which can be a non-trivial operation: consequently this approach makes easier and faster even the offline evaluation of the problem, presented in the following Section.


\subsection{Offline step}\label{subsec:Offline}
\begin{figure}
 	\centering
 	\includegraphics[width=0.6\textwidth]{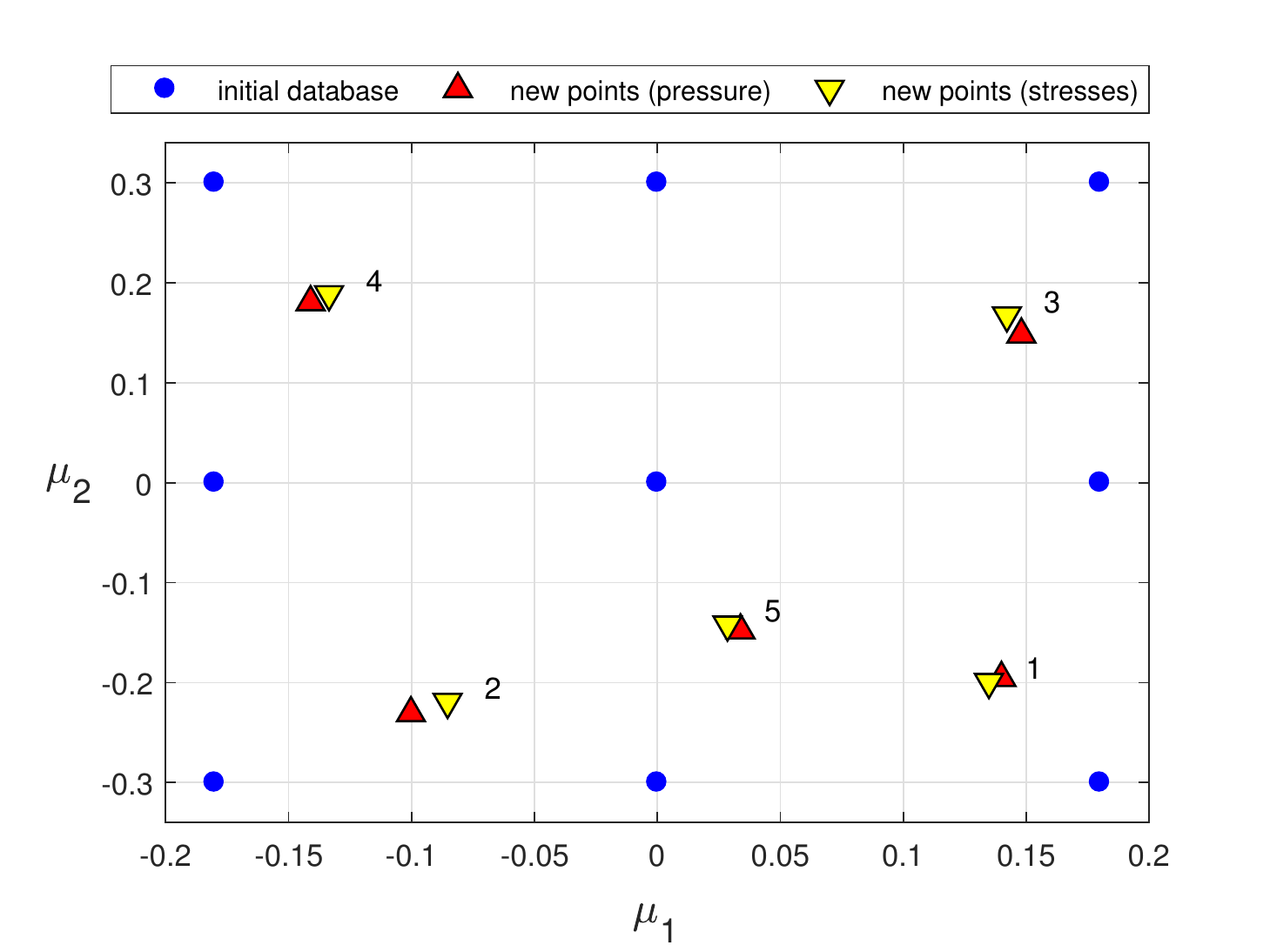}
 	\caption{Initial points and new parameter values for each iterations of the sampling algorithm.}
 	\label{fig:points}
\end{figure}

\begin{figure}
	\centering
	\subfigure[Iteration 0]{
		\setcounter{subfigure}{1}
		\includegraphics[width=0.4\textwidth]{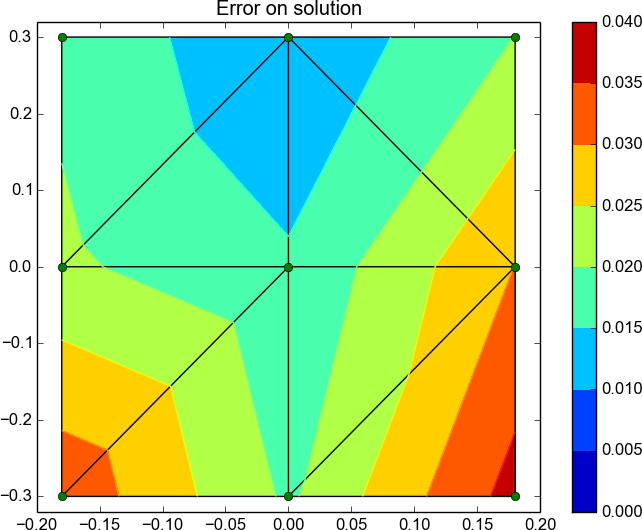}} \hspace{0.5cm}
	\subfigure[Iteration 4]{
		\setcounter{subfigure}{2}
		\includegraphics[width=0.4\textwidth]{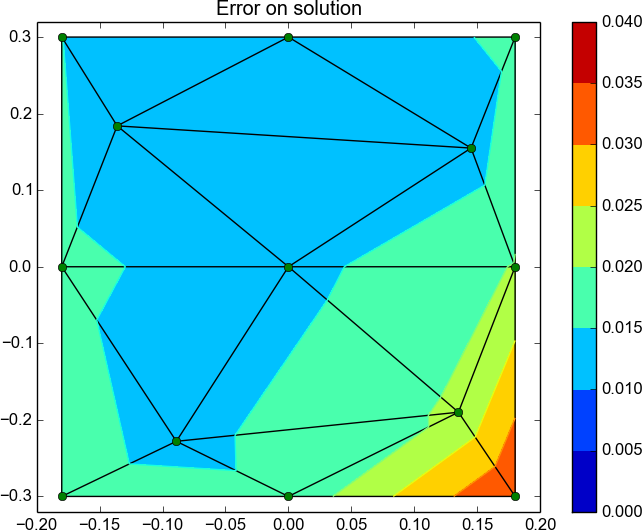}}
	\caption{L2 relative error in the parameter space for the pressure field.}
	\label{fig:Error sampling pressure}
\end{figure}

\begin{figure}
	\centering
	\subfigure[Iteration 0]{
		\setcounter{subfigure}{1}
		\includegraphics[width=0.4\textwidth]{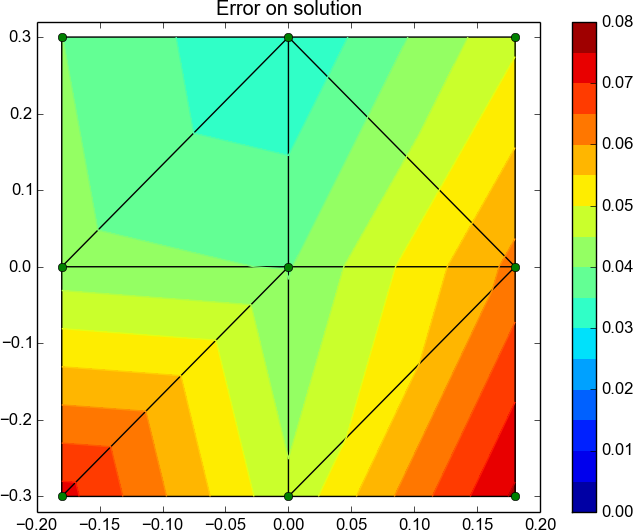}} \hspace{0.5cm}
	\subfigure[Iteration 4]{
		\setcounter{subfigure}{2}
		\includegraphics[width=0.4\textwidth]{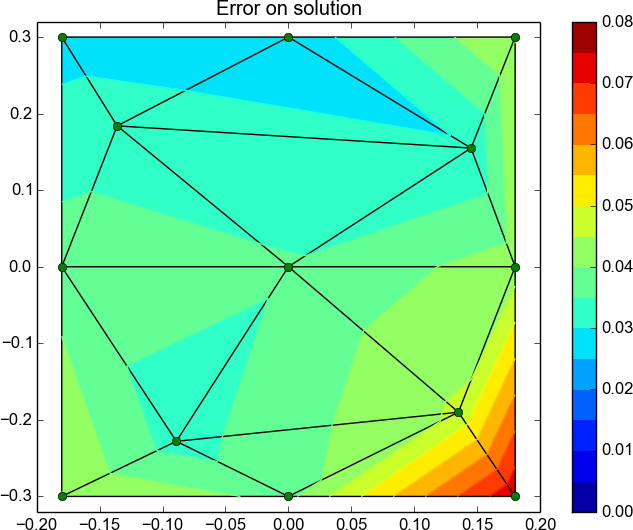}} 
	\caption{L2 relative error in the parameter space for the stresses field.}
	\label{fig:Error sampling stresses}
\end{figure}

In order to initialize both the algorithms (PODI and DD-POD), we need to evaluate a database of solutions.
As shown in Figure~\ref{fig:points}, we start from a set of $9$ points uniformly distributed in the parameter space and we iteratively add more points as chosen by Algorithm~\ref{alg:leave1out}.
Since we are interested in $p_\mathrm{w}$ and $\bm{\tau}_\mathrm{w}$ on the car in order to evaluate the drag coefficient, we compute the error estimator $e_s$ using only the restriction of the solution to the surface of the model.  
The results of the parameter selection are reported in Figures~\ref{fig:Error sampling pressure}-\ref{fig:Error sampling stresses}, where we plot the initial and the final error estimator all over the parametric domain.
It is worth noting that the algorithm chooses new parameter values which are very close using an error estimator build either on $p_\mathrm{w}$ or on $\bm{\tau}_\mathrm{w}$.
This means the physics in a given parametric point can be approximated by the physics in the already computed points in the same way for pressure and stresses.
Thus, we decided to evaluate the solution only in the point identified by the pressure error estimator: it is slightly suboptimal but allow us to halve the number of simulations.
In Table~\ref{tab:iterations}, we report the maximum and average error for each iteration: the most critical zone is the south-east quadrant (see Figures~\ref{fig:Error sampling pressure} and \ref{fig:Error sampling stresses}) but the algorithm tends to explore also regions with fewer points, resulting in average improving of the outputs reconstruction all over the parametric domain.  
\begin{table}
	\centering
	\tbl{Error estimator $\max e_s$, average error $\bar{e}_s$ and new parameter values $\mu$ for each iteration of the algorithm}
	{\begin{tabular}{r r r c r r c}
		\toprule
		& \multicolumn{3}{c}{Pressure} & \multicolumn{3}{c}{Stresses} \\
		\colrule
		Iteration & $\max e_s$ & $\bar{e}_s$ & $\bm{\mu}$ & $\max e_s$ & $\bar{e}_s$  & $\bm{\mu}$ \\
		\colrule
		0 & 0.03692 & 0.02329 & $[ 0.1399, -0.1953]$ & 0.07592 & 0.05056 & $[ 0.1348, -0.1994]$\\
		1 & 0.03697 & 0.02255 & $[-0.1002, -0.2305]$ & 0.07528 & 0.05036 & $[-0.0853, -0.2193]$\\
		2 & 0.03544 & 0.01932 & $[ 0.1480,  0.1483]$ & 0.07411 & 0.04474 & $[ 0.1421,  0.1673]$\\
		3 & 0.03535 & 0.01789 & $[-0.1410,  0.1806]$ & 0.07411 & 0.04149 & $[-0.1335,  0.1885]$\\
		4 & 0.03508 & 0.01684 & $[ 0.0338, -0.1477]$ & 0.07408 & 0.03947 & $[ 0.0285, -0.1419]$\\
		\botrule
	\end{tabular}}
	\label{tab:iterations}
\end{table}

\begin{figure}
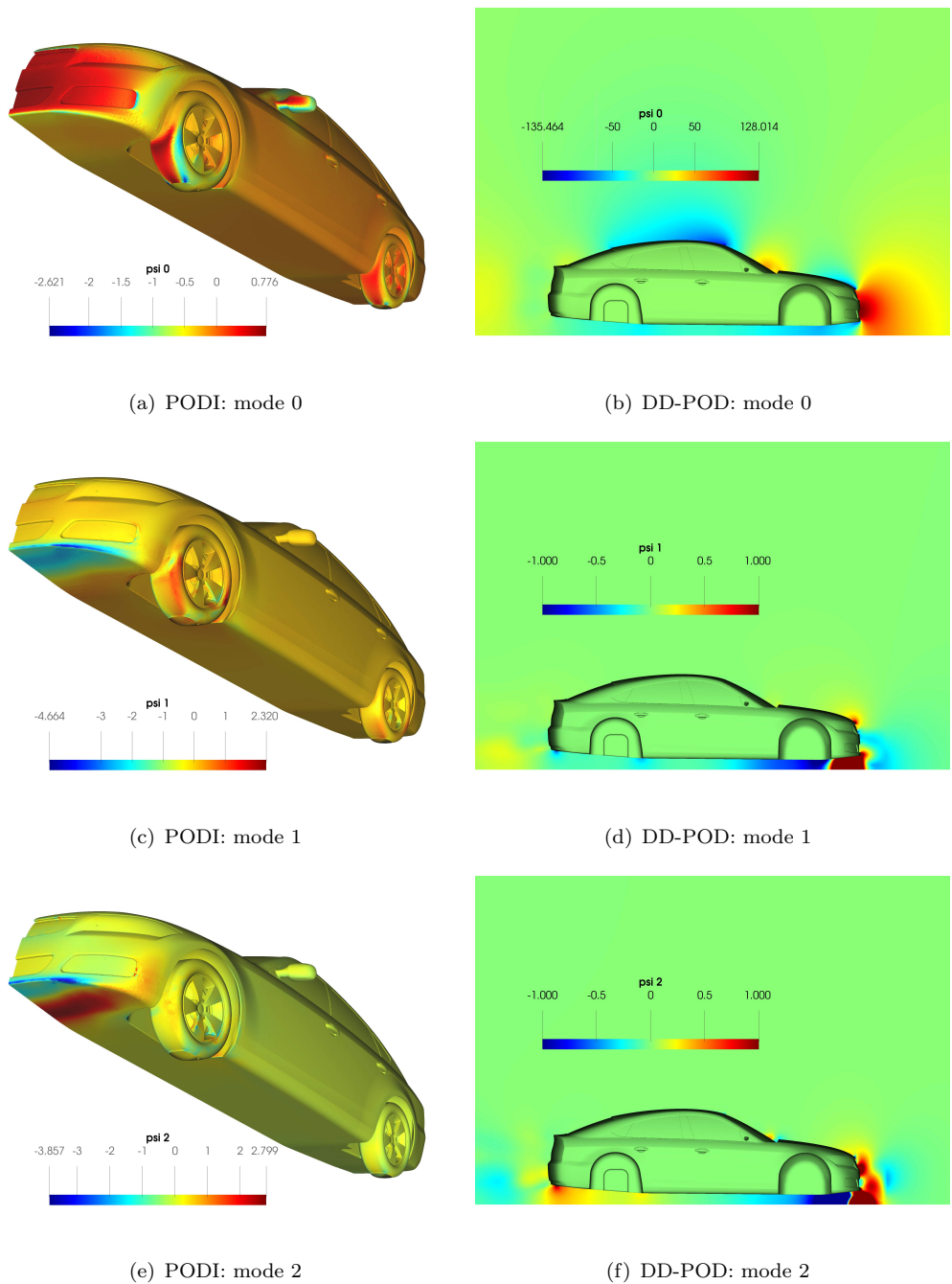

	\centering
	\subfigure[PODI: mode 0]{
		\setcounter{subfigure}{1}
		\includegraphics[width=0.4\textwidth]{figures/pod0-eps-converted-to.pdf}} 
	\subfigure[DD-POD: mode 0]{
		\setcounter{subfigure}{2}
		\includegraphics[width=0.4\textwidth]{figures/dpod0-eps-converted-to.pdf}}
	\subfigure[PODI: mode 1]{
		\setcounter{subfigure}{3}
		\includegraphics[width=0.4\textwidth]{figures/pod1-eps-converted-to.pdf}}
	\subfigure[DD-POD: mode 1]{
		\setcounter{subfigure}{4}
		\includegraphics[width=0.4\textwidth]{figures/dpod1-eps-converted-to.pdf}}
	\subfigure[PODI: mode 2]{
		\setcounter{subfigure}{5}
		\includegraphics[width=0.4\textwidth]{figures/pod2-eps-converted-to.pdf}}
	\subfigure[DD-POD: mode 2]{
		\setcounter{subfigure}{6}
		\includegraphics[width=0.4\textwidth]{figures/dpod2-eps-converted-to.pdf}}
	\caption{First three POD modes for the pressure: PODI (a,c,e) and DD-POD (b,d,f).}
	\label{fig:Pressure pod modes}
\end{figure}

\begin{figure}
	\centering
	\includegraphics[width=0.9\textwidth]{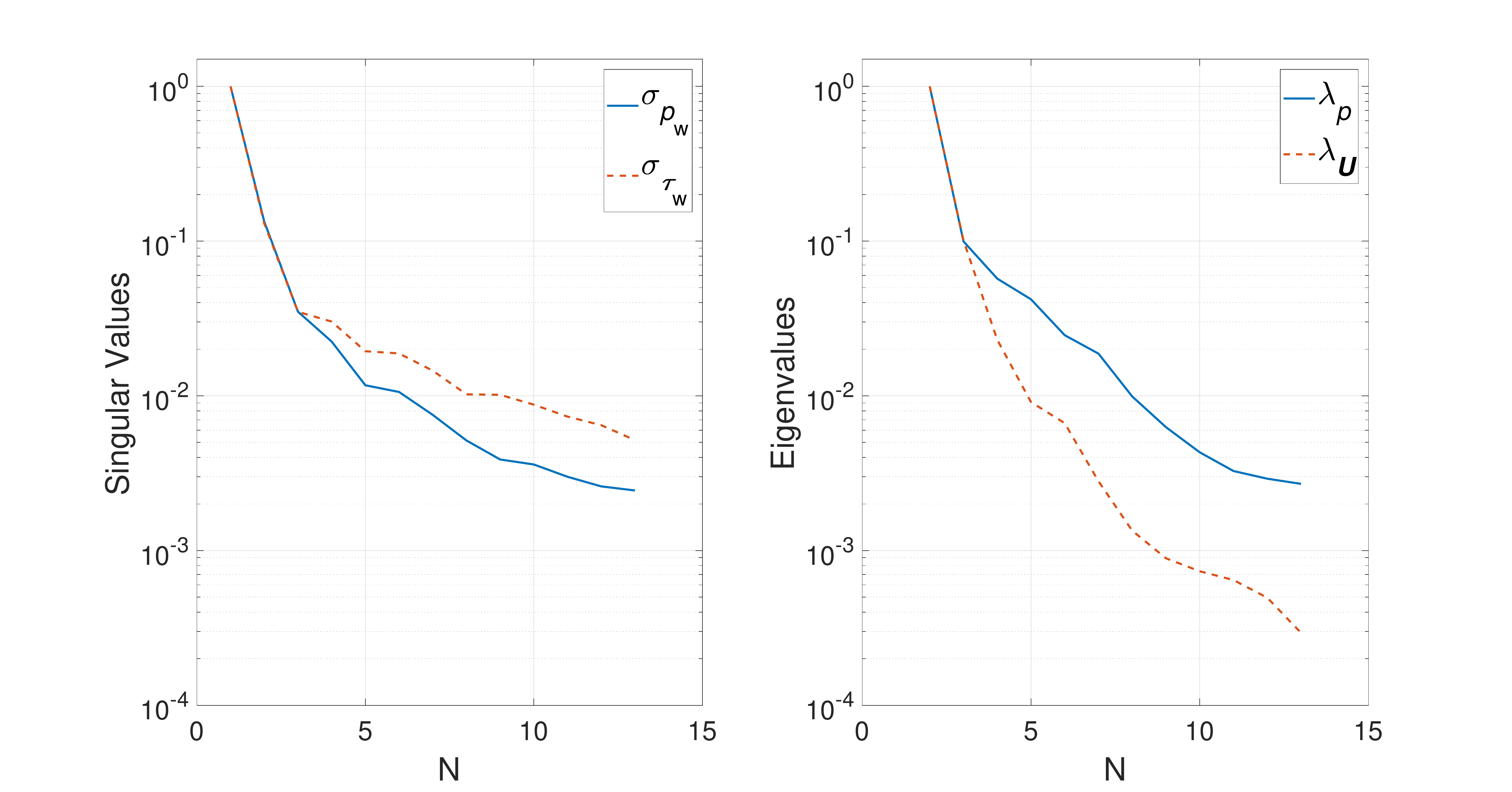}
	\caption{POD Singular Values for $p_\mathrm{w}$ and $\bm{\tau}_\mathrm{w}$ (left, PODI approach) vs POD eigenvalues for $p$ and $\bm{U}$ (right, DD-POD approach).}
	\label{fig:singular values}
\end{figure}

Once all the 13 solutions are computed, we extract the POD bases we will use in the online step and we build the triangulation of the POD coefficients over the parametric domain (for the PODI strategy).
For the PODI strategy, we build the reduced-order model on the solution only in the region of interest (i.e. the car surface), in order to speed-up the online evaluation, whereas for the DD-POD approach the POD modes are defined in the whole domain (see Figure~\ref{fig:Pressure pod modes}).
In PODI, the modes are calculated only for $p_\mathrm{w}$ and $\bm{\tau}_\mathrm{w}$; in DD-POD, since we use the POD reconstruction to impose the boundary conditions, a separated basis is computed for $p$, $\bm{U}$ and Spalart-Allmaras turbulent quantity $\tilde{\nu}$.
The decomposition for the DD-POD is chosen according to the criterion proposed by \cite{scardigli2015enabling}, resulting in a domain $\Omega_1$ of about $2\cdot10^6$ cells which includes the deformable part of the mesh.
In figure~\ref{fig:singular values} we plot both the singular values and the POD eigenvalues for the fields of interest.

\subsection{Online step} \label{subsec:Online}
In order to validate the reliability of the procedure, we evaluate four configurations out of database, chosen according to a near-random criterion (i.e. latin hypercube sampling).
Such configurations are: $\bm{\mu}=[0.135,-0.19]$, $\bm{\mu}=[-0.089,-0.228]$, $\bm{\mu}=[0.145,0.155]$ and $\bm{\mu}=[-0.136,0.184]$.

\begin{figure}
	\centering
	\includegraphics[width=0.9\textwidth]{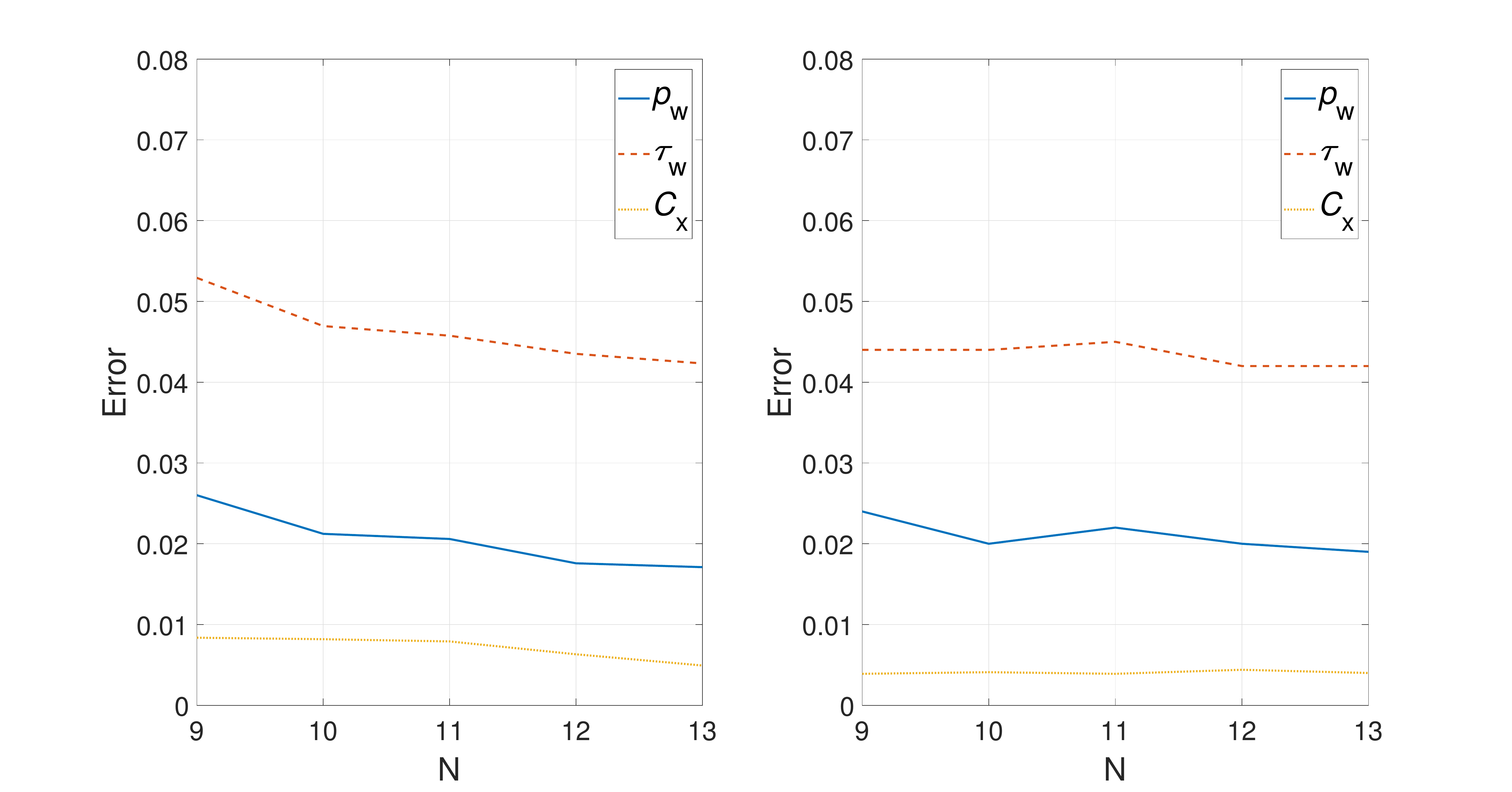}
	\caption{Error convergence enriching the POD bases for the outputs of interest: PODI (left) and DD-POD (right).}
	\label{fig:errorConvergence}
\end{figure}

In Figure~\ref{fig:errorConvergence} we plot the convergence of the averaged errors according to the dimension of the database and, subsequently, the number of POD basis.
The errors are computed as $\parallel u-u^N \parallel/ \parallel u^\mathrm{base}\parallel$ for the fields of interest and as $|C_\mathrm{x}-C_\mathrm{x}^N| / |C_\mathrm{x}^\mathrm{base}|$  for the drag coefficient, normalized with respect to the baseline configuration (i.e. $\bm{\mu}=[0,0]$).
The errors are bigger for the reconstruction of the fields, in particular the vectorial ones, as $\bm{U}$ and $\bm{\tau}_\mathrm{w}$.
This behaviour is intuitive, since in the vectorial case reduced order models have to reconstruct not only the magnitude but also the direction.
On the other hand, the error is smaller for the drag coefficient, meaning that the integration is compensating the errors, leading to a better approximation.

\begin{figure}
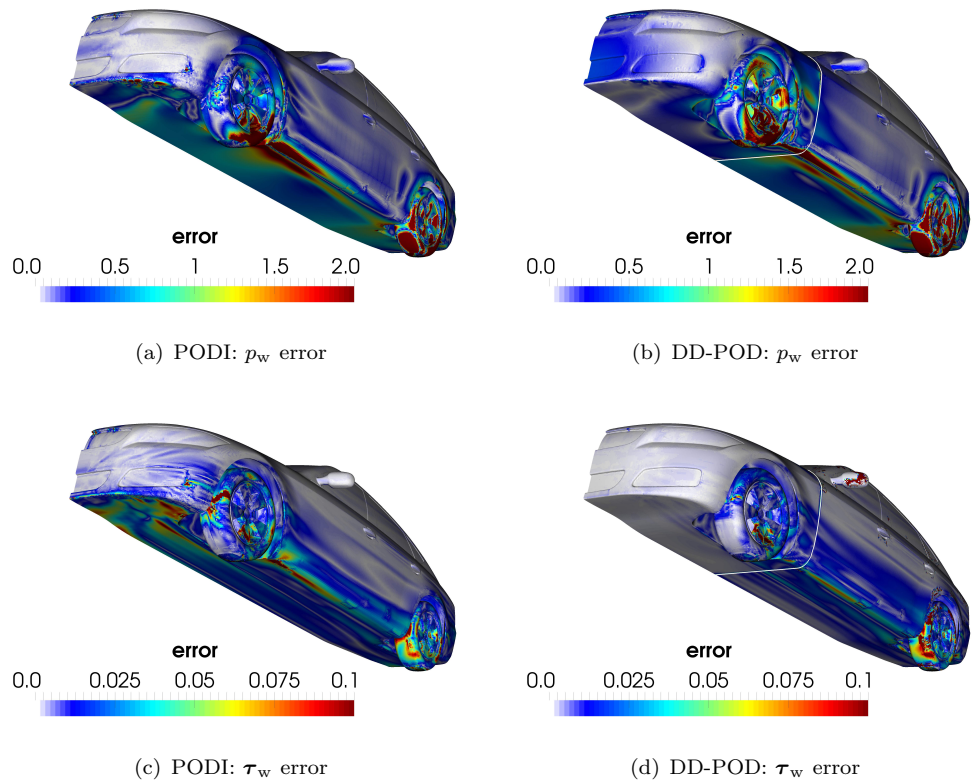

	\centering
	\subfigure[PODI: $p_\mathrm{w}$ error]{
		\setcounter{subfigure}{1}
		\includegraphics[width=0.4\textwidth]{figures/i_perr-eps-converted-to.pdf}} 
		\subfigure[DD-POD: $p_\mathrm{w}$ error]{
		\setcounter{subfigure}{2}
		\includegraphics[width=0.4\textwidth]{figures/dd_perr-eps-converted-to.pdf}} 
	\subfigure[PODI: $\bm{\tau}_\mathrm{w}$ error]{
		\setcounter{subfigure}{3}
		\includegraphics[width=0.4\textwidth]{figures/i_serr-eps-converted-to.pdf}} 
		\subfigure[DD-POD: $\bm{\tau}_\mathrm{w}$ error]{
		\setcounter{subfigure}{4}
		\includegraphics[width=0.4\textwidth]{figures/dd_serr-eps-converted-to.pdf}} 
	\caption{$p_\mathrm{w}$ and $\bm{\tau}_\mathrm{w}$ errors for the worst (in terms of accuracy) out-of-database configuration. For the DD-POD approach (right) the white line identifies the boundary between the FOM and ROM regions.}
	\label{fig:errorfields}
\end{figure}

In Figure~\ref{fig:errorfields} we plot the error between the full-order and reduced-order solution for the worst out-of-database configuration.

\begin{table}
	\centering
	\tbl{Results of online step of PODI and DD-POD framework using 13 basis functions}
	{\begin{tabular}{l l r r}
			\toprule
			& & PODI & DD-POD \\
			\colrule
			errors & $p_\mathrm{w}$ & 1.7\% & 1.9\%   \\		
			  & $\bm{\tau}_\mathrm{w}$ & 4.2\% & 4.2\%  \\
			  & $C_\mathrm{x}$ & 0.5\% & 0.4\%  \\
			\colrule
			\multicolumn{2}{l}{evaluation cost}  & $O(10^{-2} s)$ & $O(30 \, cpuh)$ \\
			\multicolumn{2}{l}{computational speed-up} & $O(10^7)$ & $O(10)$ \\
			\botrule
	\end{tabular}}
	\label{tab:performance}
\end{table}


\section{Conclusions and future work}\label{sec:conclusion}
The present work have shown a new complete pipeline for the model order reduction.
This reduction operates directly on the output of interest.
First, we have introduced a strategy for the geometry morphing of parametrized shapes, based on direct mesh deformation, showing the strength and weakness of the approach.
Then, we have introduced two different strategies for the model order reduction, both based on a POD techniques but differing in the exploitation of the extracted POD modes.
Finally we have tested and validate the present framework with an industrial benchmark coming from the automotive field, that is, the DrivAer model. 
In particular, in view of automatic shape optimization, we have been able to reconstruct the pressure and shear stress field on the surface of the model and drag coefficient with acceptable errors and a considerable speed-up, with respect to the high-fidelity solver, i.e. \textit{OpenFOAM}.
Nevertheless, we highlight that both the strategies proposed do not rely on the chosen discretization method: in particular the PODI approach treats the high-fidelity solver completely as a black box, whereas the DD-POD one requires only runtime access to the boundary conditions. 
This feature allows to exploit the user preferred software, even commercial ones.
Thus, this new framework could easily be integrated in every technical design pipeline with not much further effort.
Following this rationale, in the next future, we want to apply the present pipeline to others industrial fields, such as, for instance, in naval, nautical and aerospace engineering. 
These fields, in fact, even if they can appear very different, present very similar features: high Reynolds numbers, complex geometries, rotating mechanical parts (propellers and wheels) and the very same output of interest (pressure and wall shear stresses).
The main open issue that we have to tackle and improve is the choice of the starting sampling: the dimension of the initial database can become quickly prohibitive if we want to use a large number of design parameters.
Still, the present tool is ready-to-use for several industrial and biomedical applications.

\section*{Acknowledgements}
This work was supported by the UBE -- Underwater Blue Efficiency project (PAR-FSC programme Regione Friuli Venezia Giulia); INDAM-GNCS; and Automobili Lamborghini S.p.A..

\bibliographystyle{gCFD}
\bibliography{bibliography}

\end{document}